\documentclass[twocolumn]{autart}  
\usepackage{graphics, euscript, float,blkarray, listings}
\usepackage {graphicx}
\usepackage{tabularx, booktabs, array, float, framed}
\usepackage{color} 
\usepackage{appendix, setspace}
\usepackage{mathtools}
\usepackage[T1]{fontenc}
\usepackage[utf8]{inputenc}
\definecolor{mygreen}{RGB}{28,172,0} 
\definecolor{mylilas}{RGB}{170,55,241}
\usepackage{amsmath,amsfonts,amssymb}
\usepackage{mathtools,thmtools,etoolbox}
\usepackage{algorithm,algorithmicx,algpseudocode}
\usepackage{booktabs,subcaption}
\usepackage{tikz,pgfplots}
\pgfplotsset{compat=newest}
\usepackage{lineno}
\usepackage{todonotes}
\usepackage[colorlinks,citecolor=black]{hyperref}
\usepackage[scientific-notation=true]{siunitx}
\usepackage{multicol,cite}
\usepackage[font={large}]{caption}
\usepackage{subcaption}
\usepackage{breqn}

\newtheorem{theorem}{Theorem}[section]

\newtheorem{corollary}[theorem]{Corollary}

\theoremstyle{remark}
\theoremstyle{definition}
\newtheorem{remark}[theorem]{Remark}
\AtEndEnvironment{remark}{\hfill\ensuremath{\diamondsuit}}
\AtEndEnvironment{example}{\hfill\ensuremath{\bigcirc}}

\newcommand{\bfphi}{{\boldsymbol{\phi}}}
\newcommand{\bflambda}{{\boldsymbol{\lambda}}}

\newcommand{\bfb}{{\boldsymbol{b}}}
\newcommand{\bfc}{{\boldsymbol{c}}}

\newcommand{\goyal}{GR}

\def\mathcal{\EuScript}

\numberwithin{equation}{section}
\newtheorem{Theorem}{Theorem}[section]
\newtheorem{Lemma}[Theorem]{Lemma}
\newtheorem{Definition}[Theorem]{Definition}

\newcommand{\FH}{FHIRKA}
\newcommand{\R}{\mathbb{R}}

\newcommand{\CC}{\mathbb{C}}
\newcommand{\x}{\mathbf{x}}
\newcommand{\A}{\mathbf{A}}

\newcommand{\Y}{\mathbf{Y}}

\newcommand{\B}{\mathbf{B}}
\newcommand{\I}{\mathbf{I}}
\newcommand{\G}{\mathbf{G}}
\newcommand{\C}{\mathbf{C}}

\newcommand{\M}{\mathbf{M}}

\newcommand{\U}{\mathbf{U}}

\newcommand{\y}{\mathbf{y}}

\newcommand{\uu}{\mathbf{u}}
\newcommand{\HH}{\mathbf{H}}

\newcommand{\J}{\mathcal{J}}
\newcommand{\h}{\mathbf{h}}
\newcommand{\g}{\mathbf{g}}
\newcommand{\V}{\mathbf{V}}
\newcommand{\W}{\mathbf{W}}
\newcommand{\tf}{t_f}

\newcommand{\rc}{\boldsymbol{r}}
\newcommand{\lc}{\boldsymbol{\ell}}
\newcommand{\cc}{\boldsymbol{c}}
\newcommand{\bb}{\boldsymbol{b}}
\newcommand{\ch}{\mathcal{H}}

\newcommand{\norm}[1]{\left\lVert#1\right\rVert}
\newcolumntype{x}[1]{>{\centering\arraybackslash\hspace{0pt}}p{#1}}
\newcommand{\pder}[2][]{\dfrac{\partial#1}{\partial#2}}

\begin{document}

\begin{frontmatter}
\title{ {$\ch_2(\tf)$ Optimality Conditions for a Finite-time Horizon\thanksref{mytitlenote}} }
\author{Klajdi Sinani\thanksref{cor}} \ead{klajdi@vt.edu},
\author{Serkan~Gugercin} \ead{gugercin@vt.edu}

\thanks[mytitlenote]{This work was funded by the U.S. National Science Foundation under grants DMS-1522616 and  DMS-1720257.}
\thanks[cor]{Corresponding author }
\address{Department of Mathematics, Virginia Polytechnic Institute and State University, Blacksburg, VA 24061, USA}

\begin{abstract} 
In this paper we establish the interpolatory model reduction framework for optimal approximation of MIMO dynamical systems with respect to the $\ch_2$ norm over a finite-time horizon, denoted as the $\ch_2(\tf)$ norm.
Using the underlying inner product space, we derive the interpolatory first-order necessary optimality conditions for approximation in the $\ch_2(\tf)$ norm. Then,  we  develop an algorithm, which yields a locally optimal reduced model that satisfies the established interpolation-based optimality conditions. We test the algorithm on various numerical examples to illustrate its performance.
\end{abstract}


\end{frontmatter}
\section{Introduction}
Simulation, design, and control of dynamical systems play an important role in numerous scientific and industrial tasks such as \textcolor{black}{signal propagation in the nervous system\cite{kellems2009low};  the synthesis  of
interconnect \cite{BonD07} and semiconductor devices \cite{Hess2014}; large-scale inverse problems,
 \cite{Druskin2011solution,Lieberman2010,de2015nonlinear}; and prediction of major weather events \cite{Ant05}.}
The need for  detailed  models due to the increasing demand for greater resolution leads to  
large-scale dynamical systems, posing tremendous computational difficulties when applied in numerical simulations. In order to overcome these challenges, we perform  model reduction where we replace the large-scale dynamics with high-fidelity reduced representations.

Consider the linear time-invariant dynamical system:
\begin{equation}  \label{fom1}
\begin{aligned}
\dot{\x}(t) & = \A\x(t) +\B\uu(t),~~~\x(0) = \mathbf{0},\\
\y(t)&= \C\x(t)  = \int_{0}^{t}\h (t-\tau)\uu(\tau) d\tau,
\end{aligned}
\end{equation}
\noindent
where 
$ \A \in \R^{n\times n}$,
 $\B\in \R^{n\times m}$, and $\C\in \R^{p\times n}$ are constant matrices; 
the variable $\x(t)\in \R^n$ denotes the internal variables, $\uu(t)\in \R^m$ denotes the control inputs, and $\y(t)\in \R^p$ denotes the outputs; and $\h(t)=\C e^{\A t}\B$ is the impulse response of the full model.   The length of the internal variable $\x(t)$, i.e., $n$, is called the order of the full model that we would like to reduce.  Model
reduction achieves this by replacing the original model with a lower dimensional one:
\begin{equation}  \label{rom1}
\begin{aligned}
\dot{\x}_r(t)&=\A_r\x_r(t)+\B_r\uu(t),~~~\x_r(0) = \mathbf{0},\\
\y_r(t)&=\C_r\x_r(t)  = \int_{0}^{t}\h_r (t-\tau)\uu(\tau)d\tau,
\end{aligned}
\end{equation}
where as in \eqref{fom1}, 
$\h_r(t)=\C_r e^{\A_r t}\B_r$  is the impulse response of the reduced model, and $\A_r \in \R^{r\times r}$, $\B_r\in \R^{r\times m}$, and $\C_r\in \R^{p\times r}$ with $r\ll n$. The goal is that the output of the reduced model, $\y_r(t)$, approximates the true output, $\y(t)$, of the original system accurately in an appropriate norm. 

For the linear dynamical systems we consider here, a plethora of methods exists
for producing high-fidelity/optimal reduced models, such as  balanced truncation \cite{MulR76,Moo81} and its variants, optimal Hankel norm approximation \cite{Glo84}, and the Iterative Rational Krylov Algorithm (IRKA) \cite{GugBA08} and its variants; see \cite{Ant05,BauBF14} for further references.
These methods usually focus on constructing high-quality reduced models  over an infinite time horizon. However, in various settings, we might either have access to simulations over a finite horizon or can only simulate the system under investigation for a finite horizon such as in the case of unstable dynamical systems. Therefore, in those situations we are interested in the behavior of the dynamical system over a finite time interval $[0,t_f]$ where $t_f < \infty$, and we need the reduced model to be accurate only in the interval of interest.

Time-limited balanced truncation \cite{GawJ90,GugA03,RedK17,Kur17} and Proper Orthogonal Decomposition (POD)   \cite{HolLB96}  are two common frameworks to create reduced models on a finite horizon. 
For time-limited balanced truncation, \cite{GugA03} establishes an upper bound for the $\ch_\infty$ error between the full and reduced models,  \cite{RedK17} provides an $\ch_2$ error bound.

  
\textcolor{black}{In this paper,  we explore optimal model reduction over a finite time horizon. We use a time-limited version of the $\ch_2$ norm, denoted by $\ch_2(\tf)$, to quantify the model reduction error.}
 Optimality requires a parametrization of the reduced model. We will work with the time-domain representation of the dynamical system to derive the optimality conditions. Specifically, we represent the impulse response of the reduced dynamical system using the modal decomposition, i.e.,
 \begin{equation}\label{kernel}
\begin{aligned}
\h_r(t)=\C_r e^{\A_r t}\B_r=\sum_{i=1}^r e^{\lambda_it}\lc_i\rc_i^T.
\end{aligned}
\end{equation}
where $\lambda_i$'s are the eigenvalues of $\A_r$, and $\lc_i \in \CC^{p\times 1}, \rc_i\in \CC^{m\times1}$. 
In other words, the impulse response is expressed as a sum of $r$ rank-1 $p\times m$  matrices.
To simplify the presentation,  we assume  that $\lambda_i$'s, the reduced order poles,  are simple.  The representation \eqref{kernel} 
is nothing but a state-space transformation on $\h_r(t)=\C_r e^{\A_r t}\B_r$ using the eigenvectors of $\A_r$. 
Using the parametrization of the reduced model in \eqref{kernel}, we derive interpolatory optimality conditions in the $\ch_2(\tf)$ norm and implement a model reduction algorithm that satisfies these optimality conditions.

\textcolor{black}{The advantage of the interpolation framework we will develop  is that 
we do not require the reduced-model to be obtained via projection, as usually assumed in model reduction \cite{Ant05}. Indeed, one observation we will make is that unlike in the infinite-horizon $\ch_2$ approximation problem, the optimal reduced model in the finite-horizon case is not necessarily given by a projection and thus a projection-based approach will not be able to satisfy the optimality conditions. Therefore,  
by treating the poles and residues in \eqref{kernel} as the parameters and directly working with them, we  obtain a reduced model to satisfy the optimality conditions exactly.}


The rest of the paper is organized as follows: In Section  \ref{sec:intro} we briefly review optimal  $\ch_2$ model reduction in the infinite horizon case. The main results, including the new optimality conditions for finite horizon, are established in Section \ref{sec:main} followed by numerical examples in Section \ref{sec:num}. The papers ends with conclusions and future work in Section \ref{sec:conc}.

\section{$\ch_2$-Optimal Model Reduction: The Infinite Horizon Case} \label{sec:intro}
Model reduction with respect to the $\ch_2$ norm in the infinite horizon case has been studied extensively; see, for examples, \cite{BarCO91,BryC90,FulO90,MeiL67,Hal92,HylB85,SpaMM92,YanL99,LepMPV91,GugBA08,AniBGA13,Wil70,vuillemin2014poles,CasL18,panzer2013,gerstner2007hom,VanGA08,breiten2013near} and the references therein. 
In this section, we briefly recall these results as they will help to highlight the similarities to and differences from the finite-horizon case that we are interested in.

\subsection{$\ch_2$ Norm and $\ch_2$ Error Measure}
The error analysis for model reduction of linear dynamical systems can be conducted either in the frequency domain or in the time domain. Therefore, we define the $\ch_2$ norm in each domain. 
\begin{Definition}
Let $\h(t)$ and $\g(t)$ be the impulse responses of two asymptotically stable\footnote{We will call
$\h(t) = \C e^{\A t} \B$ asymptotically stable if all the eigenvalues of $\A$ have negative real parts. We will call 
$\h(t)$ stable when $\A$ has some semi-simple eigenvalues on the imaginary axis in addition to those with negative real parts. Otherwise, we call  $\h(t)$ unstable.}
linear dynamical systems with real state-space realizations.
The $\ch_2$ inner product $\langle \cdot, \cdot \rangle_{\ch_2}$ and the $\ch_2$ norm $\| \cdot \|_{\ch_2}$  are 
\begin{align*}
\langle \h, \g \rangle_{\ch_2} &=\int_0^\infty \mathsf{Tr}\left((\h(t))^T \g(t)\right) dt,\\
\norm{\h}_{\ch_2}&= \sqrt{ \int_0^\infty \norm{\h}_F^2 dt },
\end{align*} 
respectively, where $ \mathsf{Tr}(\cdot)$ denotes the trace and $\norm{\cdot}_F$ denotes the Frobenius norm of a matrix.
\end{Definition}
To define the $\ch_2$ norm in the frequency domain, let $\Y(s)$, $\U(s)$, and $\HH(s)$ denote the Laplace transforms of 
 the output $\y(t)$, the input $\uu(t)$, and the impulse response $\h(t) = \C e^{\A t}\B$
in \eqref{fom1}. Then, taking the Laplace transform of the convolution integral in  \eqref{fom1}, we obtain

$\Y(s)= \HH(s) \U(s),  \quad \mbox{where} \quad \HH(s) =\C(s\I-\A)^{-1}\B$

is called the transfer function of \eqref{fom1}. Let $\{\rho_1,\rho_2,\ldots,\rho_n\}$ denote the eigenvalues of $\A$, assumed simple. Then, similar to the parametrization of the reduced model $\h_r(t)$ in \eqref{kernel}, the full-model impulse response can be equivalently written as 
\begin{equation} \label{eqn:hfull}
\h(t) = \sum_{i=1}^n e^{\rho_i t} \cc_i \bb_i^T~~\mbox{with}~~
\HH(s) =\sum_{i=1}^n \frac{\cc_i \bb_i^T}{s-\rho_i},
\end{equation}
where $\cc_i \in \CC^{p}$ and $\bb_i \in \CC^m$, for $i=1,\ldots,n$.
This is called the pole-residue form where $\rho_i$'s are the poles of the (rational) transfer  function 
$\HH(s)$ with the corresponding rank-$1$ residues $\cc_i \bb_i^T$.
\begin{Definition}
Let $\HH(s)$ and $\G(s)$ denote the transfer functions  of two asymptotically stable  dynamical systems with real state-space realizations.
The $\ch_2$ inner product $\langle \cdot, \cdot \rangle_{\ch_2}$ 
and the $\ch_2$ norm $\| \cdot \|_{\ch_2}$ are
\begin{align*}
\langle \G, \HH \rangle_{\ch_2}&:=\dfrac{1}{2\pi}\int_{-\infty}^\infty \mathsf{Tr}(\overline{\G}(-\imath\omega)\HH^T(\imath\omega)) d\omega \\
\norm{ \HH}_{\ch_2}&:=\sqrt{\dfrac{1}{2\pi} \int_{-\infty}^\infty \norm{\HH(i\omega)}_F^2 d\omega},
\end{align*}
respectively.
\end{Definition}

Similar to $\HH(s)$, let $\HH_r(s)=\C_r(s\I_r-\A_r)^{-1}\B_r$ denote the transfer function of the reduced model. Then, the relevance and importance of the $\ch_2$ norm in the model reduction become clear by noting that 
\begin{equation}
\begin{aligned}
\norm{\y-\y_r}_{L_\infty}\leq\norm{\HH-\HH_r}_{\ch_2}\norm{\uu}_{L_2},
\end{aligned}
\end{equation}
i.e., the $L_\infty$ norm of the output error $\y(t)-\y_r(t)$ due to model reduction is bounded by the 
the $\ch_2$ norm of error transfer function relative to the $L_2$ norm of the input $\uu(t)$;
see, e.g.,  \cite{AntBG10b}, for a proof.
Therefore, to guarantee that the reduced model output $\y_r(t)$ is close to the original one $\y(t)$, one might look for a reduced model that minimizes the $\ch_2$ error norm.

\subsection{Interpolatory Conditions for Optimal $\ch_2$ Model Reduction}
Given a reduced order $r$, the goal is to construct a reduced order model whose transfer function $\HH_r(s)$ minimizes the $\ch_2$ error norm $\norm{\HH - \HH_r}_{\ch_2}$.  Since this is a non-convex optimization problem, the usual, numerically feasible, approach  is to find a reduced model that satisfies the necessary conditions for $\ch_2$ optimality. These conditions can be formulated in terms of Sylvester equations \cite{Wil70,HylB85} or interpolation \cite{MeiL67,GugBA08}. These two frameworks are equivalent \cite{GugBA08}. In this paper, we will focus on the interpolation framework.
\begin{theorem} \label{h2thm}
Let 
$$
\h_r(t)=\sum_{k=1}^r e^{\lambda_kt}\lc_k\rc_k^T \qquad \Longleftrightarrow \qquad \HH_r(s) = \sum_{k=1}^r \frac{\lc_k\rc_k^T}{s-\lambda_k}
$$
be the best $r^{th}$ order approximation of an asymptotically stable linear dynamical system  $\HH(s)$ with respect to the $\ch_2$ norm. Then, for $k=1, 2, ... , r$,
\begin{equation} \label{eqn:h2cond}
\begin{aligned}
\lc_k^T\HH(-\lambda_k)&=\lc_k^T\HH_r(-\lambda_k),\\~~
\HH(-\lambda_k)\rc_k&=\HH_r(-\lambda_k)\rc_k,\\
\lc_k^T\HH'(-\lambda_k)\rc_k&=\lc_k^T\HH_r'(-\lambda_k)\rc_k,
\end{aligned}
\end{equation}

{\color{black} where $\HH'(s)$ denotes the derivative of $\HH(s)$ with respect to $s$.}
\end{theorem}
{\color{black}For more details on Theorem \ref{h2thm}, see {\cite {GugBA08,AntBG10b}}.}
This result states that the transfer function of the optimal $\ch_2$ approximation to $\HH(s)$ is a (tangential) Hermite interpolant where the interpolation points are the mirror images of the reduced-order poles $\{\lambda_k\}$, and the tangental directions are given by the corresponding residues $\{\lc_k\rc_k^T\}$.  Since the optimality conditions depend on the reduced-model to be computed, the solution requires a nonlinear iteration. The Iterative Rational Krylov Algorithm (IRKA)  \cite{GugBA08} and its variants such as \cite{beattie2012realization,gerstner2007hom,CasL18,beattie2009trust} use these interpolation based optimality conditions to produce   an interpolatory, locally $\ch_2$ optimal reduced model. The next section will extend this framework to the finite-time interval case.

\section{$\ch_2(\tf)$ Optimal Model Reduction on a Finite Horizon}  \label{sec:main} 
%
{\color{black}
In this section, we present  the main theoretical results of the paper i.e., the interpolatory $\ch_2(\tf)$  optimality conditions, and discuss their implications.}


\subsection{$\ch_2(t_f)$ Norm on a Finite-time Horizon}
It is immediately clear from the time-domain definition of the (infinite-horizon) $\ch_2$ norm how to define the finite-horizon version:
\begin{Definition}\label{innerh2}
Let $\h(t)$ and $\g(t)$ denote the impulse responses of two dynamical systems with real state-space realizations.  
For a finite-time horizon 
$[0,t_f]$, the $\ch_2(\tf)$\footnote{\textcolor{black}{Even though the term  $\ch_2$ is mostly associated with a measure in the frequency domain, following \cite{RedK17,GoyR17} we are using the notation $\ch_2(\tf)$ here as well to denote the error measure specifically formulated in the time domain. Our main reason is to keep the connection to the infinite horizon problem where the measure in the frequency- and time-domains are equivalent. And more importantly, as in the regular $\ch_2$ case, the optimality conditions will still appear as interpolation conditions in the frequency domain. 
 }}
   inner product $\langle \cdot, \cdot \rangle_{\ch_2(\tf)}$ 
and $\ch_2(\tf)$ norm $\norm{\cdot}_{\ch_2(\tf)}$ are defined as  
\begin{align*}
\langle \h, \g \rangle_{\ch_2(\tf)} &=\int_0^{\tf} \mathsf{Tr}( (\h(t))^T \g(t)) dt, \\
\norm{\h}_{\ch_2(\tf)} &=\sqrt{\int_0^{\tf} \norm{\h(t)}_F^2 dt}.
\end{align*}

\end{Definition}

\subsection{Finite Horizon Interpolation-based Conditions $\ch_2(\tf)$ Optimal Model Reduction}
The problem we are interested in is as follows: Given a dynamical system with impulse response $\h(t)$ (or equivalently with transfer function $\HH(s)$) and a reduced order $r$, find the reduced model with the impulse response 
\begin{equation}\label{kernel2}
\begin{aligned}
\h_r(t)=\sum_{i=1}^r e^{\lambda_it}\lc_i\rc_i^T.
\end{aligned}
\end{equation}
such that $\norm{\h - \h_r}_{\ch_2(\tf)}$ is minimized. As in the regular $\ch_2$ case, this is a non-convex optimization problem and we will focus on local minimizers. Using the parametrization \eqref{kernel2}, we we will derive interpolation-based necessary conditions for optimality.
The main result is given by Theorem \ref{mh2opt}. However, we  need many supplementary results, Lemmas 
\ref{minnerexp}-\ref{mgexp}, to reach this final conclusion. It is immediately clear that, the $\ch_2(t_f)$-error, denoted by $\J$, satisfies 
\begin{equation}\label{merexp}
\begin{aligned}
\J&=\norm{\h-\h_r}_{\ch_2({\tf})}^2\\
&=\norm{\h}_{\ch_2({\tf})}^2-2 \langle\h, \h_r\rangle_{\ch_2({\tf})}+\norm{\h_r}_{\ch_2({\tf})}^2,
\end{aligned}
\end{equation}
where the inner product  $\langle\h, \h_r\rangle_{\ch_2({\tf})}$ is real since $\h(t)$ and $\h_r(t)$ are real.
Finding the first-order necessary conditions for optimal ${\ch_2({\tf})}$ model reduction will require computing the gradient of the error expression \eqref{merexp} with respect to the optimization variables. 
Since the reduced model, as described by  the impulse response in $\h_r(t)$, is parametrized by the 
reduced order poles $\{\lambda_i\}$, and the residue directions $\{\lc_i\}$ and $\{\rc_i\}$, we will be computing the gradient of the error with respect to these variables. Since the first term in the error \eqref{merexp}, i.e., $\norm{\h}_{\ch_2({\tf})}$, is a constant, we will be focusing on the remaining two terms only. First, in the next lemma, we will formulate these two last terms in terms of $\{\lambda_i\}$, 
$\{\lc_i\}$ and $\{\rc_i\}$.
\begin{Lemma}\label{minnerexp}

Let  $\h (t)= \sum_{j=1}^n e^{\rho_jt} \cc_j \bb_j^T $ and $\h_r(t) = \sum_{i=1}^r e^{\lambda_it}\lc_i  \rc_i^T $ be, respectively, the impulse responses of the full and reduced models  as described in \eqref{eqn:hfull} and \eqref{kernel2}.
Then,
\begin{align} \label{eqn:hhr}
\langle\h, \h_r\rangle_{\ch_2({\tf})}=\sum_{j=1}^n\sum_{i=1}^r\lc_{i}^T\cc_j\bb_j^T\rc_{i}\dfrac{e^{(\lambda_i+\rho_j){\tf}}-1}{\lambda_i+\rho_j},
\end{align}
and
\begin{align} \label{eqn:hrhr}
\norm{\h_r}_{\ch_2({\tf})}^2= \sum_{i=1}^r\sum_{j=1}^r\lc_{i}^T\lc_{j}\rc_{j}^T\rc_{i}\dfrac{e^{(\lambda_i+\lambda_j){\tf}}-1}{\lambda_i+\lambda_j}.
\end{align}
\end{Lemma}
\begin{pf}
The both results  follow from the definition of the $\ch_2(\tf)$ inner product.  First consider
\begin{align*}
\langle\h, \h_r\rangle_{\ch_2({\tf})}&=\mathsf{Tr}\left(\int_0^{\tf}\h_r(t)^T\h(t)\, dt\right).
\end{align*}
Plug $\h (t)= \sum_{j=1}^n \cc_j \bb_j^T e^{\rho_jt}$ and $\h_r(t) = \sum_{i=1}^r \lc_i  \rc_i^T e^{\lambda_it}$
into this formula to obtain
\begin{align*}
\langle\h, \h_r\rangle_{\ch_2({\tf})}&=\mathsf{Tr} \left(\int_0^{\tf}\sum_{i=1}^r(\lc_{i}\rc_{i}^Te^{\lambda_it})^T\sum_{j=1}^n\cc_j\bb_j^Te^{\rho_jt}dt\right)\\
 &= \mathsf{Tr} \left(\sum_{i=1}^r 
\sum_{j=1}^n 
\rc_{i}\lc_{i}^T\cc_j \bb_j^T \int_0^{\tf}e^{(\lambda_i + \rho_j )t} dt\right)
\end{align*}
Computing the integral and using the fact that $\mathsf{Tr}(\A_1 \A_2) = \mathsf{Tr}(\A_2 \A_1)$  for two matrices $\A_1$ and $\A_2$ of  appropriate sizes,   we obtain
\begin{align*}
\langle\h, \h_r\rangle_{\ch_2({\tf})}&=\mathsf{Tr}\left(\sum_{j=1}^n\sum_{i=1}^r\rc_{i}\lc_{i}^T\cc_j\bb_j^T\dfrac{e^{(\lambda_i+\rho_j){\tf}}-1}{\lambda_i+\rho_j}\right)\\
&=\sum_{j=1}^n\sum_{i=1}^r\lc_{i}^T\cc_j\bb_j^T\rc_{i}\dfrac{e^{(\lambda_i+\rho_j){\tf}}-1}{\lambda_i+\rho_j},
\end{align*}
which proves \eqref{eqn:hhr}. Then, \eqref{eqn:hrhr} follows directly by replacing $\h(t)$ with $\h_r(t)$ in this derivation. 
\end{pf}
For infinite time horizon, Theorem \ref{h2thm} tells us that a locally $\ch_2$ optimal reduced transfer function is a tangential Hermite interpolant of the original transfer function at the mirror images of the reduced poles. We will show that in the finite horizon case, even though Hermite tangential interpolation is still the necessary conditions for optimality,  what is being interpolated and what the interpolant is are different.

\begin{theorem}\label{mh2opt}
Let $\HH(s)= \C(s\I-\A)^{-1} \B$, with $\A \in \R^{n \times n}$, $\B \in \R^{n \times m}$, and $\C \in \R^{p \times n}$,   be the transfer function of the full-order model with  the pole-residue representation 
$
\HH(s) =  \sum_{i=1}^n \frac{\cc_i \bb_i^T}{s-\rho_i}
$
as in \eqref{eqn:hfull}, where $\cc_i \in \CC^{p}$, $\bb_i \in \CC^m$, and $\rho_i \in \CC$ for $i=1,\ldots,n$. 
For a finite-time horizon $[0,t_f]$, define 
\begin{equation} \label{eqn:defG}
\G(s)=-e^{-s\tf}\C(s\I-\A)^{-1}e^{\A \tf}\B+\HH(s).
\end{equation}
Let 
\begin{equation}
\HH_r(s)= \C_r(s\I_r-\A_r)^{-1} \B_r = \sum_{i=1}^r \frac{\lc_i \rc_i^T}{s-\lambda_i} 
\end{equation}
be the transfer function of  the best $r^{th}$ order approximation of $\HH(s)$ with respect to the $\ch_2(\tf)$ norm
where $\A_r \in \R^{r \times r}$, $\B_r \in \R^{r \times m}$, $\C \in \R^{p \times r}$,  $\lc_i \in \CC^{p}$, $\rc_i \in \CC^m$, and $\lambda_i \in \CC$ for $i=1,\ldots,r$.
Define 
\begin{equation} \label{eqn:defGr}
\G_{r}(s)=-e^{-s\tf}\C(s\I_r-\A_r)^{-1}e^{\A_r \tf}\B_r+\HH_r(s).
 \end{equation}
 Then,  for $k=1, 2, ... , r$,
\begin{align}
\lc_k^T\G(-\lambda_k)&=\lc_k^T\G_r(-\lambda_k),  \label{eqn:cond1} \\
\G(-\lambda_k)\rc_k&=\G_r(-\lambda_k)\rc_k,~~\mbox{and}  \label{eqn:cond2} \\
\lc_k^T\G'(-\lambda_k)\rc_k&=\lc_k^T\G_r'(-\lambda_k)\rc_k. \label{eqn:cond3}
\end{align}
\end{theorem}

The next lemma will be used in proving Theorem \ref{mh2opt}.
\begin{Lemma}\label{mgexp}
Let $\G(s)$ and $\G_r(s)$ be as defined in \eqref{eqn:defG} and \eqref{eqn:defGr}, respectively. Then,
\begin{align}
\G(-\lambda_k)&=\phantom{-}\sum_{j=1}^n\cc_j\bb_j^T\dfrac{e^{(\lambda_k+\rho_j){\tf}}-1}{\lambda_k+\rho_j},
\label{eqn:Gatlk}\\
 \G'(-\lambda_k)& =-\sum_{j=1}^n\cc_i\bb_i^T\dfrac {({\tf} (\lambda_k+\rho_j)-1)e^{(\lambda_k+\rho_j){\tf}}+1}{(\lambda_k+\rho_j)^2}, \label{eqn:Gpatlk}  \\
 \G_r(-\lambda_k)&=\phantom{-}\sum_{j=1}^n\lc_j\rc_j^T\dfrac{e^{(\lambda_k+\lambda_j){\tf}}-1}{\lambda_k+\lambda_j}, ~~\mbox{and} 
\label{eqn:Gratlk}\\
 \G_r'(-\lambda_k)& =-\sum_{j=1}^n\lc_j\rc_j^T\dfrac {({\tf} (\lambda_k+\lambda_j)-1)e^{(\lambda_k+\lambda_i){\tf}}+1}{(\lambda_k+\lambda_j)^2}. \label{eqn:Grpatlk} 
\end{align}
\end{Lemma}
\begin{pf}
Recall the definition of $\G(s) = -e^{-s\tf}\C(s\I-\A)^{-1}e^{\A \tf}\B+\HH(s)$. Note that we assume that the eigenvalues of $\A$ are simple. Therefore, $e^{\A \tf}$ is also diagonalizable by the eigenvectors of $\A$. Using this fact and the pole-zero residue decomposition of   $\HH(s) = \C(s\I-\A)^{-1} \B = \sum_{j=1}^n \frac{\cc_j \bb_j^T}{s-\rho_j}$, we obtain 
\begin{align}
\G(s) &= -e^{-s\tf}\C(s\I-\A)^{-1}e^{\A \tf}\B+\HH(s)  \nonumber \\ 
&=-e^{-s{\tf}}\sum_{j=1}^n\cc_j\bb_j^T\dfrac{e^{\rho_j\tf}}{s-\rho_j}+\sum_{j=1}^n\cc_j\bb_j^T\dfrac{1}{s-\rho_j}\\
&=\sum_{j=1}^n\cc_j\bb_j^T\dfrac{-e^{(-s + \rho_j){\tf}}}{s-\rho_j}+\sum_{j=1}^n\cc_j\bb_j^T\dfrac{1}{s-\rho_j} \nonumber\\
&=\sum_{j=1}^n\cc_j\bb_j^T\dfrac{e^{(-s+\rho_j)\tf}-1}{-s+\rho_j}. \label{eq:Gs}
\end{align}
Thus,
$\G(-\lambda_k)=\sum_{j=1}^n\cc_j\bb_j^T\dfrac{e^{(\lambda_k+\rho_j){\tf}}-1}{\lambda_k+\rho_j}$,
which proves \eqref{eqn:Gatlk}. To prove \eqref{eqn:Gpatlk}, we first differentiate \eqref{eq:Gs} with respect to $s$ to obtain
 \begin{align*}
 \G'(s) =\sum_{j=1}^n\cc_j\bb_j^T\dfrac {{\tf} (s-\rho_j)e^{(-s+\rho_j){\tf}}+e^{(-s+\rho_j){\tf}}-1}{(s-\rho_j)^2}
 \end{align*}
 Plugging in $s= -\lambda_k$ in this last expression yields the desired result \eqref{eqn:Gpatlk}. The proofs of
 \eqref{eqn:Gpatlk} and \eqref{eqn:Grpatlk} follow analogously. 
\end{pf}

\textbf{{PROOF OF THEOREM \ref{mh2opt}.}}
 As mentioned above, let $\J$ denote the $\ch_2({\tf})$ error norm square, i.e., 
\begin{align*}
 \J&=\norm{\h -\h_r}^2_{\ch_2({\tf})} \\
&=\norm{\h}_{\ch_2({\tf})}^2-2 \langle\h, \h_r\rangle_{\ch_2({\tf})}+\norm{\h_r}_{\ch_2({\tf})}^2.
\end{align*}
The expressions for $ \langle\h, \h_r\rangle_{\ch_2({\tf})}$ and $\norm{\h_r}_{\ch_2({\tf})}^2$ in terms of the optimization variables $\{\lambda_k\}$, $\{\rc_k\}$, and $\{\lc_k\}$, for $k=1,2,\ldots,r$, were already derived in Lemma \eqref{minnerexp}. To make the gradient computations with respect to the $k$th parameter more clear, we seperate the $k$th term from these expressions. For example, we write $\langle\h, \h_r\rangle_{\ch_2({\tf})}$ in \eqref{eqn:hhr} as
\begin{dmath*}
\langle\h, \h_r\rangle_{\ch_2({\tf})}=\sum_{j=1}^n\lc_{k}^T\cc_j\bb_j^T\rc_{k}\dfrac{e^{(\lambda_k+\rho_j){\tf}}-1}{\lambda_k+\rho_j}
+\sum_{j=1}^n\sum_{\substack{i=1\\i\neq k}}^r\lc_{i}^T\cc_j\bb_j^T\rc_{i}\dfrac{e^{(\lambda_i+\rho_j){\tf}}-1}{\lambda_i+\rho_j}.
\end{dmath*}
Following the same procedure for $\norm{\h_r}_{\ch_2({\tf})}^2$, we obtain
\begin{dmath} \label{merror}
\J= \int_0^{\tf}\h(t)^T\h(t) dt-2 \left(\sum_{j=1}^n\lc_{k}^T\cc_j\bb_j^T\rc_{k}\dfrac{e^{(\lambda_k+\rho_j){\tf}}-1}{\lambda_k+\rho_j}
+\sum_{j=1}^n\sum_{\substack{i=1\\i\neq k}}^r\lc_{i}^T\cc_j\bb_j^T\rc_{i}\dfrac{e^{(\lambda_i+\rho_j){\tf}}-1}{\lambda_i+\rho_j}\right) 
+ \lc_{k}^T\lc_{k}\rc_{k}^T\rc_{k}\dfrac{e^{(2\lambda_k){\tf}}-1}{2\lambda_k}
+\sum_{\substack{i=1 \\i \neq k}}^r\lc_{i}^T\lc_{k}\rc_{k}^T\rc_{i}\dfrac{e^{(\lambda_i+\lambda_k){\tf}}-1}{\lambda_i+\lambda_k}
+\sum_{\substack{j=1\\j\neq k}}^r\lc_{k}^T\lc_{j}\rc_{j}^T\rc_{k}\dfrac{e^{(\lambda_k+\lambda_j){\tf}}-1}{\lambda_i+\lambda_j}+\sum_{\substack{j=1\\j\neq k}}^r\sum_{\substack{i=1\\i\neq k}}^r\lc_{i}^T\lc_{j}\rc_{j}^T\rc_{i}\dfrac{e^{(\lambda_i+\lambda_j){\tf}}-1}{\lambda_i+\lambda_j}.
\end{dmath}
To compute the gradient of the cost function $\J$ we perturb the cost functional with respect to the residue directions, i.e., $\lc_k\to \lc_k+\Delta\lc_k$ and $\rc_k\to \rc_k+\Delta\rc_k$: 
\begin{dmath*}
\widetilde{\J}_k= \int_0^{\tf}\h(t)^T\h(t) dt-2\bigg( \sum_{j=1}^n(\lc_{k}+\Delta\lc_{k})^T\cc_j\bb_j^T(\rc_{k}+\Delta\rc_{k})\dfrac{e^{(\lambda_k+\rho_j){\tf}}-1}{\lambda_k+\rho_j} 
+\sum_{j=1}^n\sum_{\substack{i=1\\i\neq k}}^r\lc_{i}^T\cc_j\bb_j^T\rc_{i}\dfrac{e^{(\lambda_i+\rho_j){\tf}}-1}{\lambda_i+\rho_j}\bigg)+ (\lc_{k}+\Delta\lc_{k})^T(\lc_{k}+\Delta\lc_{k})(\rc_{k}+\Delta\rc_{k})^T(\rc_{k}+\Delta\rc_{k})\dfrac{e^{(2\lambda_k){\tf}}-1}{2\lambda_k}
+\sum_{\substack{i=1 \\i \neq k}}^r\lc_{i}^T(\lc_{k}+\Delta\lc_{k})(\rc_{k}+\Delta\rc_{k})^T\rc_{i}\dfrac{e^{(\lambda_i+\lambda_k){\tf}}-1}{\lambda_i+\lambda_k}+                             
\sum_{\substack{j=1\\j\neq k}}^r(\lc_{k}+\Delta\lc_{k})^T\lc_{j}\rc_{j}^T(\rc_{k}+\Delta\rc_{k})\dfrac{e^{(\lambda_k+\lambda_j){\tf}}-1}{\lambda_i+\lambda_j}
+\sum_{\substack{j=1\\j\neq k}}^r\sum_{\substack{i=1\\i\neq k}}^r\lc_{i}^T\lc_{j}\rc_{j}^T\rc_{i}\dfrac{e^{(\lambda_i+\lambda_j){\tf}}-1}{\lambda_i+\lambda_j}.
\end{dmath*}
Then, collecting  the terms that are multiplied by $\Delta\lc_{k}$ and $\Delta\rc_{k}$, we obtain 
\begin{dmath*}
\nabla_{\rc_{k} }\J =-2\lc_k^T\sum_{j=1}^n\cc_j\bb_j^T\dfrac{e^{(\lambda_k+\rho_j){\tf}}-1}{\lambda_k+\rho_j}+2\lc_k^T\sum_{j=1}^n\lc_j\rc_j^T\dfrac{e^{(\lambda_k+\lambda_j){\tf}}-1}{\lambda_k+\lambda_j},
\end{dmath*}
\begin{dmath*}
\nabla_{\lc_{k} }\J=-2\left (\sum_{j=1}^n\cc_j\bb_j^T\dfrac{e^{(\lambda_k+\rho_j){\tf}}-1}{\lambda_k+\rho_j}\right)\rc_k+2\left(\sum_{j=1}^n\lc_j\rc_j^T\dfrac{e^{(\lambda_k+\lambda_j){\tf}}-1}{\lambda_k+\lambda_j}\right)\rc_k.\\
\end{dmath*}
Setting $\nabla_{\rc_{k} }\J=0   \text{ and }\nabla_{\lc_{k} }\J=0 $, and using Lemma \ref{mgexp}, mainly
\eqref{eqn:Gatlk} and \eqref{eqn:Gratlk}, yield 
\begin{align*}
\lc_k^T\G(-\lambda_k)&=\lc_k^T\G_r(-\lambda_k), \\ 
\G(-\lambda_k)\rc_k&=\G_r(-\lambda_k)\rc_k,
\end{align*}
which proves \eqref{eqn:cond1} and \eqref{eqn:cond2}. To prove \eqref{eqn:cond3},  we differentiate $\J$ in \eqref{merror}  with respect to the $k$-th pole $\lambda_k$. Note that we have written $\J$ in such a way to isolate the terms that depend on $\lambda_k$ from the ones that do not. Thus, many of the terms in \eqref{merror} have a zero derivative and we obtain:
\begin{dmath}\label{eqn:diffl}
\pder[\J]{\lambda_k}=-2\lc_k^T \left(\sum_{j=1}^n \cc_j\bb_j^T \dfrac{{\tf}(\lambda_k+\rho_j)e^{(\rho_j+\lambda_k){\tf}}} {(\lambda_k+\rho_j)^2}     -\dfrac       {e^{(\rho_j+\lambda_k){\tf}}-1}{(\lambda_k+\rho_j)^2}\right)\rc_k 
+2\lc_k^T \left(\sum^r_{i=1} \lc_i\rc_i^T\dfrac{{\tf}(\lambda_i+\lambda_j\k)e^{(\lambda_i+\lambda_k){\tf}}}  {(\lambda_i+\lambda_k)^2}         -\dfrac{e^{(\lambda_i+\lambda_k){\tf}}-1}{(\lambda_i+\lambda_k)^2}\right)\rc_k.
\end{dmath}
{\color{black}
Note that the first term in \eqref{eqn:diffl} corresponds to the derivative of the second term in \eqref{merror} and the second term in \eqref{eqn:diffl} corresponds to the derivative of the last four terms in \eqref{merror}.}
We rewrite \eqref{eqn:diffl} to obtain 

\begin{dmath}
\pder[\J]{\lambda_k} =-2K_1+2K_2
\end{dmath}
where
\begin{dmath}\label{eqn:l11} 
K_1=\lc_k^T \left(\sum_{j=1}^n \cc_j\bb_j^T \dfrac{({\tf}(\lambda_k+\rho_j)-1)e^{(\rho_j+\lambda_k){\tf}}+1}{(\lambda_k+\rho_j)^2}\right)\rc_k  
\end{dmath}
\begin{dmath}\label{eqn:l12} 
K_2=\lc_k^T \left(\sum^r_{i=1} \lc_i\rc_i^T\dfrac{({\tf}(\lambda_i+\lambda_j\k)-1)e^{(\lambda_i+\lambda_k){\tf}}+1}{(\lambda_i+\lambda_k)^2}\right)\rc_k.
\end{dmath}

Lemma \eqref{mgexp}, specifically \eqref{eqn:Gpatlk} and \eqref{eqn:Grpatlk}, show that the expressions in the parentheses in \eqref{eqn:l11}  and \eqref{eqn:l12}   are, respectively, $-\G'(-\lambda_k)$ and $-\G_r'(-\lambda_k)$.
If $\pder[\J]{\lambda_k}= 0$, then
\begin{align*}
\lc_k^T\G'(-\lambda_k)\rc_k&=\lc_k^T\G_r'(-\lambda_k)\rc_k,
\end{align*}
which completes the proof.
\begin{flushright} $\Box$  \end{flushright}
{\color{black} We note that the interval of interest is problem dependent and the choice of the interval, i.e., the choice of $t_f$, may depend on the model. However, these optimality conditions hold for any choice of $t_f>0$.}

\begin{remark}
{\color{black} In the infinite-horizon case, if  $\HH_r(s)$ is the best $\ch_2$ approximation to $\HH(s)$, then $\HH_r(s)$ interpolates $\HH(s)$. However, in the finite-horizon case, 
the interpolant is $\G_r(s)$, and the interpolated function is $\G(s)$; thus $\HH_r(s)$ does not interpolate $\HH(s)$. To give more intuition about these resulting interpolation conditions, consider the time-limited function $\g(t)$ such that 
$\g(t)=\h(t)$ when   $t\leq {\tf}$ and $\g(t)=0$ when $t > {\tf}$.
A direct calculation shows that $\G(s)$ is the Laplace transform of $\g(t)$. Similarly let $\g_r(t)$ 
denote the time-limited version of $\h_r(t)$. Then its Laplace transform is 
 $\G_r(s)$.  Therefore, the optimality conditions  \eqref{eqn:cond1}--\eqref{eqn:cond3}
 correspond to optimal interpolation of $\G(s)$ (Laplace transform of the time-limited function
 $\g(t)$)
 by $\G_r(s)$ (Laplace transform of the time limited function  $\g_r(t)$). The fact that  
 $\g(t)$ and   $\g_r(t)$ are both time-limited is the precise reason  why we cannot simply apply $\ch_2$ optimal reduction  to $\G(s)$. The method of \cite{beattie2012realization}, called TF-IRKA, does not require the original function to be a rational function. Thus, in principle we can use TF-IRKA to reduce $\G(s)$. However, the resulting reduced model is a rational function without any structure. In our case, the reduced model $\G_r(s)$ needs to retain the same structure as $\G(s)$ so that we can extract an $\HH_r(s)$. In other words, if we simply apply an $\ch_2$ optimal algorithm to $\G(s)$, we would be approximating a finite horizon model by an infinite horizon one and we cannot extract $\HH_r(s)$. Therefore, a new  algorithmic framework  is needed as we discuss in more detail in Section \ref{sec:discuss}.}
\end{remark}

\begin{remark}
{\color{black}
For an unstable dynamical system without any purely imaginary eigenvalues, one can work with the
 $\mathcal{L}_2$ norm by decomposing it into a stable and anti-stable system, and then obtain an interpolatory reduced model based on this measure.  However, this solution requires destroying the causality of the underlying dynamics \cite{magruder2010rational}. This is not the framework we are interested in here and we  work with a finite-time interval. }
 \end{remark}

\subsection{Implication of the interpolatory $\ch_2({\tf})$ optimality conditions}
\label{sec:discuss}
Theorem \ref{mh2opt} extends the interpolatory infinite-horizon $\ch_2$ optimality conditions \eqref{eqn:h2cond} to the finite-horizon case. Note that in the case of asymptotically stable dynamical systems, if we let $t_f \to \infty$, we recover the infinite-horizon conditions  \eqref{eqn:h2cond}.

 The major difference from the regular $\ch_2$ problem is that  optimality no longer requires that the reduced model $\HH_r(s)$ tangentially interpolate the full model $\HH(s)$. Instead, the auxiliary reduced-order function $\G_r(s)$ should be a tangential Hermite interpolant to the auxiliary  full-order  function $\G(s)$. However, the optimal interpolation points and the optimal tangential directions still result from the pole-residue representation of the reduced-order transfer function $\HH_r(s)$. This situation is  similar to the interpolatory optimality conditions for the  frequency-weighted $\ch_2$-optimal model reduction problem in which  one tries to minimize a weighted $\ch_2$ norm in the frequency domain, i.e.,
find $\HH_r(s)$ that minimizes $ \norm{ \mathbf{W} ( \HH - \HH_r) }_{\ch_2}$ where $\mathbf{W}(s)$ represents a weighting function in the frequency domain. As \cite{breiten2013near} showed, the optimality  
in the frequency-weighted $\ch_2$-norm requires that a function of $\HH_r(s)$ tangentially interpolate
a function of  $\HH(s)$. Despite this  conceptual similarity, the resulting interpolation conditions are drastically different from what we obtained here as one would expect due to the different measures. For details, we refer the reader to  \cite{breiten2013near}.

As we pointed out in Section \ref{sec:intro}, in addition to the interpolatory framework, one can represent the $\ch_2$ optimality conditions in terms of  Sylvester equations, leading to a projection framework for the reduced model. This means that given the full-model $\HH(s) = \C(s\I - \A)^{-1} \B$, one constructs two bases $\V,\W \in \R^{n\times r}$ with $\V^T\W = \I_r$ such that the reduced-order quantities are obtained via projection, i.e.,
\begin{equation} \label{eqn:proj}
\A_r = \W^T \A \V,~~\B_r = \W^T \B, \quad \mbox{and}\quad \C_r = \C \V.
\end{equation}
In the infinite-horizon case, \cite{Wil70} showed that the optimal $\ch_2$ reduced model is indeed guaranteed to be obtained via projection. Recently, Goyal and Redmann  \cite{GoyR17} have established the Sylvester-equation based optimality conditions for the time-limited $\ch_2$ model reduction problem; i.e., they extended the Wilson framework to the time-limited (finite-horizon) $\ch_2$ problem. Furthermore, they have developed a projection-based IRKA-type numerical algorithm to construct the reduced models. However, as the authors point out in \cite{GoyR17}, even though their algorithm yields high-fidelity reduced models in terms of the $\ch_2(t_f)$ measure, the resulting reduced model  satisfies the optimality conditions only approximately. This is not surprising in light of the optimality conditions we derived here. Since the optimality requires that $\G_r(s)$ should interpolate $\G(s)$ (as opposed to $\HH_r(s)$ interpolating $\HH(s)$), unlike in the infinite-horizon case, the reduced model in the finite-horizon case is \emph{not} necessarily given by a projection  as in \eqref{eqn:proj}. Therefore, a projection-based approach would satisfy the optimality conditions only approximately. This was also the case in  \cite{breiten2013near} where even though a projection-based IRKA-like algorithm produced high-fidelity reduced models in the weighted norm,  it satisfied the optimality conditions approximately. 

The advantage of the interpolation framework and the parametrization \eqref{kernel} we consider here is that
we do not require the reduced-model to be obtained via projection. By treating the poles and residues in \eqref{kernel} as the parameters and directly working with them, we can obtain a reduced model to satisfy the optimality conditions exactly. Even though the main focus of this paper is the theoretical interpolatory framework  and a robust numerical algorithm will be  fully considered in a future work, in the next section we will  discuss a basic numerical framework one can develop using the interpolatory conditions.

\begin{remark}
The finite-horizon approximation problem for  \emph{discrete-time dynamical systems} has been considered 
in \cite{MelVG14}. The derivation in \cite{MelVG14}, however, allows the reduced-model quantities to vary at every time- step, thus using a time-varying reduced model as opposed to the time-invariant formulation considered here and in  \cite{GoyR17}. Allowing time-varying quantities drastically simplifies the gradient computations, leading to a recurrence relations for optimality. Therefore, the model reduction problem for finite-horizon $\ch_2$ approximation for \emph{time-invariant} discrete-time dynamical systems is still an open question.
\end{remark}
\section{Numerical considerations}  \label{sec:num}
In this section, we briefly discuss a numerical framework to construct a reduced model that satisfies the optimality conditions \eqref{eqn:cond1}-\eqref{eqn:cond3}. To make the presentation and discussion concise,  we will focus on the single-input/single-output (SISO) version only. The complete numerical framework for the general case  with further details on the underlying optimization schemes will be discussed in a separate work.
\subsection{A descent-type algorithm for the single-input/single-output case}
Let $\HH(s)$ and $\HH_r(s)$ be SISO full- and reduced-model transfer functions, respectively, i.e.,
\begin{equation} \label{eqn:sisoh}
\begin{aligned}
\HH(s) &=  \bfc^T(s\I-\A)^{-1}\bfb = \sum_{i=1}^n \frac{\psi_i}{s-\rho_i}\\
\HH_r(s) &=  \bfc_r^T(s\I_r-\A_r)^{-1}\bfb_r = \sum_{i=1}^r \frac{\phi_i}{s-\lambda_i},
\end{aligned}
\end{equation}
where $\A \in \R^{n\times n}$, $\bfb,\bfc \in \R^n$, $\A_r \in \R^{r\times r}$, and 
$\bfb_r,\bfc_r \in \R^r$. Note that the residues $\psi_i$ and $\phi_i$ are scalar valued.
The following result, which is an immediate consequence of Theorem \ref{mh2opt}, summarizes the optimality conditions for SISO systems.
\begin{corollary}\label{h2opt}
Given the SISO transfer functions $\HH(s)$ and $\HH_r(s)$ as defined in \eqref{eqn:sisoh}, define

\begin{align*}
\G(s)&=-e^{-s{\tf}}\bfc^T(s\I-\A)^{-1}e^{\A {\tf}}\bfb+\HH(s),\\
\G_r(s)&=-e^{-s{\tf}}\bfc^T_r(s\I_r-\A_r)^{-1}e^{\A_r {\tf}}\bfb_r+\HH_r(s).
\end{align*}
If $\HH_r$ is the best $r^{th}$ order approximation of $\HH$ with respect to the $\ch_2({\tf})$ norm, then
\begin{align} \label{def:GGrsiso}
\G(-\lambda_k)=\G_r(-\lambda_k),\ \mbox{and}\
\G'(-\lambda_k)=\G_r'(-\lambda_k)
\end{align}
where  $\lambda_k$ for $k=1, 2, ... , r$ are the poles of the reduced system $\HH_r(s)$ as given in
 \eqref{eqn:sisoh}.
\end{corollary}
As stated before, the $\ch_2(t_f)$ minimization problem is a non-convex optimization problem and 
Corollary \ref{h2opt} gives the necessary conditions for optimality when both  poles and residues are treated as variables. However, if the poles are fixed, we can establish the necessary and sufficient optimality conditions for the residues and find the global minimizer,  the optimal residues, by solving a linear system.

\begin{corollary}\label{optres}
Let $\HH(s)$ and $\HH_r(s)$ be as given in \eqref{eqn:sisoh}, and $\G(s)$ and $\G_r(s)$ as
in \eqref{def:GGrsiso}.  Assume  the reduced  poles   $\{\lambda_i\}_{i=1}^r$ are fixed. Then,
$\HH_r(s)$ is the best $r^{th}$ order approximation of $\HH(s)$ with respect to the $\ch_2({\tf})$ norm if and only if  $\M \bfphi= \mathbf{z}$, \mbox{or~equivalently},$\G(-\lambda_k)=\G_r(-\lambda_k)$,\ \mbox{for}\ $k = 1,2,\ldots, r,$
where  $\bfphi = [\phi_1~\phi_2~\cdots~\phi_r]^T \in \CC^r$ is the vector of  residues; 
$\mathbf{z} \in \CC^r$ is the vector with  entries 
$$
\mathbf{z}_j=e^{\lambda_j{\tf}}\bfc^T(-\lambda_j\I-\A)^{-1}e^{\A {\tf}}\bfb-\HH(-\lambda_j), 
i=1,2,\ldots,r ;
$$
and 
$\M \in \CC^{r\times r}$  is the  matrix with  entries
$$
\M_{i,j}=\dfrac{e^{(\lambda_i+\lambda_j){\tf}}-1}{\lambda_i+\lambda_j}, \quad \mbox{for}\quad
i,j=1,2,\ldots,r.$$
\end{corollary}
\begin{pf}
Using the SISO counterparts for \eqref{eqn:hrhr} and \eqref{eqn:hhr} and applying some algebraic manipulation, the cost functional $\J$ can be written as
\begin{align*}
\J&= \norm{\h}_{\ch_2(t_f)}^2  - 2\bfphi^{T}\mathbf{w} + \bfphi^{T}\M\bfphi,
\end{align*}
where $\mathbf{w} \in \CC^{r\times 1}$ has the entries
 $$\mathbf{w}_i=\sum_{k=1}^n \psi_k \dfrac{e^{(\rho_k+\lambda_i){\tf}}-1}{\lambda_i+\rho_k}  \quad \mbox{for}\quad i=1,2,\ldots,r.$$
Note that 
$\M$ is positive definite, $\phi^{T}\M\phi=\norm{\h_r}_{\ch_2({\tf})}^2 > 0$, and the cost function is quadratic in $\bfphi$. Thus the (global) minimizer is  obtained by solving $\M \bfphi= \mathbf{z}$, which corresponds to rewriting  $\G(-\lambda_k)=\G_r(-\lambda_k)$ for $k = 1,2,\ldots, r$ in a compact way. 
\end{pf}
The result is analogous to the regular infinite-horizon $\ch_2$ problem where the Lagrange optimality becomes necessary and sufficient once the poles are fixed \cite{gaier1987lectures,beattie2012realization}. What is important here is that once the reduced poles are fixed, the best residues can be computed directly by solving an $r\times r$ linear system $\M \bfphi= \mathbf{z}$.  This is the property that we will exploit in the numerical scheme next.

\subsubsection{FHIRKA: A numerical algorithm for $\ch_2(t_f)$ model reduction}
Here we describe a numerical algorithm which produces a reduced model that satisfies the necessary $\ch_2({\tf})$ optimality conditions upon convergence. Let $\bflambda \in \C^r$ denote the vector of reduced poles. Thus,
the error $\J$ is a function of $\bflambda$ and $\bfphi$. Since we explicitly know the gradients of the cost function
with respect to $\bflambda$ and $\bfphi$ (and indeed we can compute the Hessians as well), one can (locally) minimize $\J$ using well established optimization tools. However, as   Corollary \ref{optres} shows,  we can easily compute the globally optimal $\bfphi$ for fixed  $\bflambda$. Therefore, 
we will treat the reduced  poles $\bflambda$ as the optimization parameter, and once $\bflambda$ are updated at the $kth$ step of an optimization algorithm, we find/update the  corresponding optimal residues $\bfphi$ based on Corollary \ref{optres}, and then repeat the process. Similar strategies have been successfully employed in the regular 
$\ch_2$ optimal approximation problem as well; see \cite{beattie2012realization,beattie2009trust}. In summary,
we use a quasi-Newton type optimization as $\bflambda$ being the parameter and in each optimization  step,
we update the residues,  $\bfphi$, by solving the $r\times r$ linear system $\M \bfphi = \mathbf{z}$  as in Corollary \ref{optres}.  Since we are enforcing interpolation at every step of the algorithm, yet tackling the model reduction problem over a finite horizon, we name this algorithm Finite Horizon IRKA, denoted by \FH.  Unlike regular IRKA,  \FH \ is a descent algorithm, thus indeed mimics \cite{beattie2009trust} more closely. Upon convergence, the locally optimal reduced model satisfies the first-order necessary conditions of 
Corollary \ref{def:GGrsiso}.

\subsection{Numerical Results}  
In this section we compare the proposed algorithm FHIRKA with  Proper Orthogonal Decomposition (POD), Time-Limited Balanced Truncation (TLBT), and the recently introduced $\ch_2(t_f)$-based algorithm  by Goyal and Redmann (\goyal) \cite{GoyR17}, as briefly discussed in Section \ref{sec:discuss}.  

We use three models: a heat model of order $n=197$ \cite{Niconet}, a model of  the International Space Station 1R Module (ISS 1R) of order $n=270$ \cite{gugercin2001approximation}, and a toy unstable model of order $n=402$. 
The ISS 1R model has $3$-inputs and $3$-outputs. We  focus on the SISO subsystem from the first-input to the first-output.
We have created the unstable system such that it has $400$ stable poles and $2$ unstable poles (positive real part).

\begin{figure}[H] 
 \centering
 \includegraphics [scale=0.35]{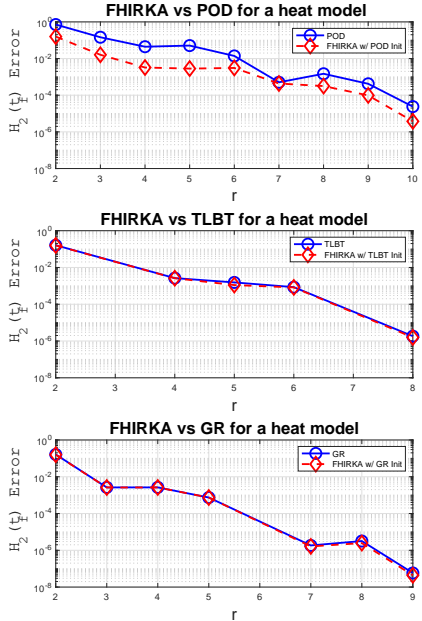} 
 \caption{\FH \ and other algorithms for the heat model\label{fig:hm}}
 \end{figure}

\begin{figure}[H]
 \centering
  \includegraphics [scale=0.35]{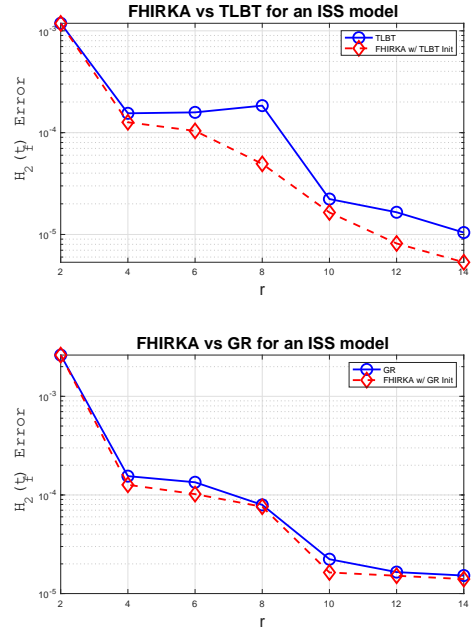}
 \caption{\FH \ and other algorithms for the ISS model \label{fig:iss}}
 \end{figure}

\begin{figure}[H]
 \centering
   \includegraphics [scale=0.35]{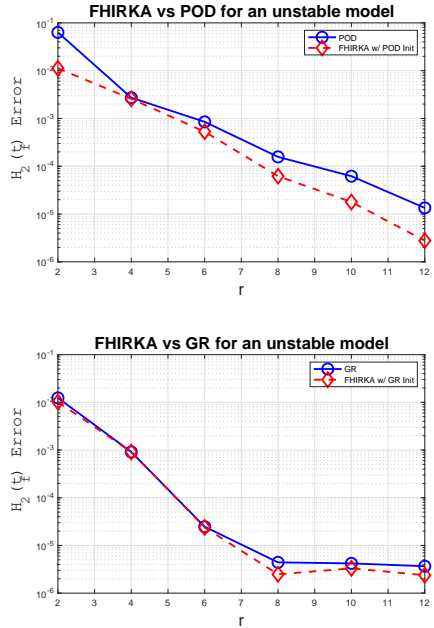}
      \caption{\FH \ and other algorithms for the unstable model \label{fig:uns}}
 \end{figure}

For all three models, we choose $t_f = 1$, first reduce the original model using POD, \goyal \  or TLBT, and then use the resulting reduced model  to initialize \FH. Thus, we are trying to investigate how these different initializations affect the final reduced model via  \FH \ and how much improvement one might expect.
The results are shown in Figures \ref{fig:hm}--\ref{fig:uns},  where we show the $\ch_2(\tf)$ approximation error for different values of $r$, the order of the reduced model.
  All three initializations are used for the heat model (Figure \ref{fig:hm}) where the order is reduced from $r=2$ to $r=10$ with increments of one. For some $r$ values, certain initializations are excluded (e.g., the \goyal \ initialization for $r=6$) since the algorithm either did not converge or produced  poor approximations. However, this happened only rarely.
  
   For the ISS model (Figure \ref{fig:iss}), we use  TLBT and \goyal \ initializations since  POD approximations was very poor and  is excluded. In this case, we reduce the order from $r=2$ to $r=14$ with increments of $2$. For the unstable model (Figure \ref{fig:uns}), we use POD and \goyal \ initializations; for this model our implementation of TLBT produced  poor results and is avoided. In this case,  we reduce the order from $r=2$ to $r=12$ with increments of $2$.  
 The first observation is that,  since  \FH \ is a descent-method and drives the initialization to a local minimizer,  it improves the accuracy of  the reduced model for all three initializations as expected. The improvements could be dramatic. For example, 
 \FH \ is able to outperform POD as much as by an order of magnitude, see, for example,
 Figure \ref{fig:hm}, the  $r=4$ and $r=5$ cases. 
 While \FH \ improves TLBT and \goyal\ initialization as well, the improvements for the heat model are not as significant. However, for the ISS model,
\FH \ is able to improve the TLBT performance as much as $50\%$; see, e.g., Figure \ref{fig:iss}, the $r=8$ case. 
The best improvement of the \goyal \ initialization has occurred for the unstable model. For example, for $r=8$, for the unstable model, \FH\ improved the reduced model by more than $40 \%$. Gains were significantly better for POD, especially for larger $r$ values. 
{\color{black} 
}


%

 
\begin{figure}[H]
 \centering
   \includegraphics [scale=0.35]{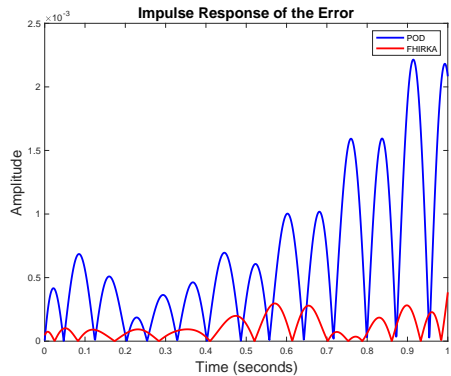}
      \caption{\FH \ and POD for the ISS model \label{fig:impulse}}
 \end{figure}

\textcolor{black}{Finally, in Figure \ref{fig:impulse}, we compare the error in the impulse responses due to POD and \FH \ for the ISS model. For both methods, POD and \FH\  the reduced model was of order $r=14$. As we can see from the graph, \FH \ clearly outperforms POD on the time interval $[0, 1]$.}

  Overall,  as expected, \FH \ yields a better approximation compared to the other algorithms for each model. 
 We find that \goyal\  provided the best initialization for \FH. This is not surprising since \goyal\ produces a reduced-model that approximately satisfies the $\ch_2(t_f)$ optimality conditions. 
 
 As we stated above, the numerical issues will be fully investigated  in a future work where we will develop a robust interpolatory $\ch_2(t_f)$-descent algorithm for MIMO systems. We will not only study better initialization techniques in terms of speed and accuracy, but also make the algorithm numerically more efficient by using approximation  techniques for the matrix exponential $e^{\A t_f}$ appearing in the $\ch_2(t_f)$ setting. 
We will also investigate the MIMO version of Corollary \ref{optres}. In the MIMO case, even for fixed poles, one cannot simply find the globally optimal residue directions by solving a linear system, since the problem is no longer quadratic in these variables. In the regular $\ch_2$ case, finding the optimal residue directions for given poles required solving a nonlinear least-squares problem \cite{beattie2012realization}. We anticipate a similar formulation here and will investigate the corresponding numerical implications.

\section{Conclusions and Future Work} \label{sec:conc}
We  established interpolatory $\ch_2({\tf})$-optimality conditions for model reduction of MIMO dynamical systems over a finite horizon. Even though the optimal interpolation points and tangential directions are still determined by the reduced model, we showed that unlike the regular $\ch_2$ problem, a modified reduced-transfer function should interpolate a modified full-order transfer function. For the special case of SISO models, we have studied a numerical algorithm and illustrated that it performs effectively.

As in the regular $\ch_2$-problem, establishing equivalency between the Sylvester-equation based $\ch_2(t_f)$-optimality conditions of \cite{GoyR17} and the interpolation-based conditions we developed here will be an interesting direction to pursue. Furthermore, extensions to bilinear and quadratic-bilinear problems will be crucial. 


\bibliographystyle{plain}

\bibliography{bibliography}

\begin{thebibliography}{10}

\bibitem{AniBGA13}
B.~Ani\`{c}, C.~Beattie, S.~Gugercin, and A.~C. Antoulas.
\newblock Interpolatory weighted-$\mathcal{H}_2$ model reduction.
\newblock {\em Automatica}, 49:1275--1280, 2013.

\bibitem{Ant05}
A.~C. Antoulas.
\newblock {\em Approximation of large-scale dynamical systems (Advances in
  Design and Control)}.
\newblock Society for Industrial and Applied Mathematics, Philadelphia, PA,
  USA, 2005.

\bibitem{AntBG10b}
A.C. Antoulas, C.~Beattie, and S.~Gugercin.
\newblock Interpolatory model reduction of large-scale dynamical systems.
\newblock In J.~Mohammadpour and K.~Grigoriadis, editors, {\em Efficient
  Modeling and Control of Large-Scale Systems}. Springer-Verlag, 2010.

\bibitem{BarCO91}
L.~Baratchart, M.~Cardelli, and M.~Olivi.
\newblock Identification and rational $\ell_2$ approximation: A gradient
  algorithm.
\newblock {\em Automat.}, 27:413--418, 1991.

\bibitem{BauBF14}
U.~Baur, P.~Benner, and L.~Feng.
\newblock Model order reduction for linear and nonlinear systems: a
  system-theoretic perspective.
\newblock {\em Archives of Computational Methods in Engineering},
  21(4):331--358, 2014.

\bibitem{beattie2009trust}
C.A. Beattie and S.~Gugercin.
\newblock A trust region method for optimal $\mathcal{H}_2$ model reduction.
\newblock In {\em Decision and Control, 2009 held jointly with the 2009 28th
  Chinese Control Conference. CDC/CCC 2009. Proceedings of the 48th IEEE
  Conference on}, pages 5370--5375. IEEE, 2009.

\bibitem{beattie2012realization}
C.A. Beattie and S.~Gugercin.
\newblock Realization-independent $\mathcal{H}_2$-approximation.
\newblock In {\em Proceedings of 51st IEEE Conference on Decision and Control},
  pages 4953 -- 4958, 2012.

\bibitem{BonD07}
B.N. Bond and L.~Daniel.
\newblock A piecewise-linear moment-matching approach to parameterized
  model-order reduction for highly nonlinear systems.
\newblock {\em IEEE Transactions on Computer-Aided Design of Integrated
  Circuits and Systems}, 26(12):2116--2129, 2007.

\bibitem{breiten2013near}
T.~Breiten, C.A. Beattie, and S.~Gugercin.
\newblock Near-optimal frequency-weighted interpolatory model reduction.
\newblock {\em Systems \& Control Letters}, 78:8--18, 2015.

\bibitem{BryC90}
A.~E. Bryson and A.~Carrier.
\newblock Second-order algorithm for optimal model order reduction.
\newblock {\em J. Guidance Control Dynam.}, 13:887--892, 1990.

\bibitem{gerstner2007hom}
A.~Bunse-Gerstner, D.~Kubali{\'n}ska, G.~Vossen, and D.~Wilczek.
\newblock $\mathcal{H}_2$-optimal model reduction for large scale discrete
  dynamical mimo systems.
\newblock {\em Journal of Computational and Applied Mathematics},
  233(5):1202--1216, 2010.

\bibitem{CasL18}
A.~Castagnotto and A.~Lohmann.
\newblock A new framework for $h_2$-optimal model reduction.
\newblock {\em Mathematical and Computer Modelling of Dynamical Systems}, pages
  1--22, 2018.

\bibitem{Niconet}
Y.~Chahlaoui and P.~Van Dooren.
\newblock A collection of benchmark examples for model reduction of linear time
  invariant dynamical systems.
\newblock Technical report, SLICOT Working Note 2002-2, 2002.

\bibitem{de2015nonlinear}
E.~De~Sturler, S.~Gugercin, M.~Kilmer, S.~Chaturantabut, C.A. Beattie, and
  M.~O'Connell.
\newblock Nonlinear parametric inversion using interpolatory model reduction.
\newblock {\em SIAM Journal on Scientific Computing}, 37(3):B495--B517, 2015.

\bibitem{Druskin2011solution}
V.~Druskin, V.~Simoncini, and M.~Zaslavsky.
\newblock Solution of the time-domain inverse resistivity problem in the model
  reduction framework part {I}. {One}-dimensional problem with {SISO} data.
\newblock {\em SIAM Journal on Scientific Computing}, 35(3):A1621--A1640, 2013.

\bibitem{FulO90}
P.~Fulcheri and M.~Olivi.
\newblock Matrix rational $\mathcal{H}_2$ approximation: A gradient algorithm
  based on schur analysis.
\newblock {\em SIAM J. Control Optim.}, 36:2103--2127, 1998.

\bibitem{gaier1987lectures}
D.~Gaier.
\newblock {\em Lectures on Complex Approximation}.
\newblock Birkhauser, Cambridge, MA, 1987.

\bibitem{GawJ90}
W.~Gawronski and J.-N. Juang.
\newblock Model reduction in limited time and frequency intervals.
\newblock {\em International Journal of Systems Science}, 21(2):349--376, 1990.

\bibitem{Glo84}
K.~Glover.
\newblock All optimal hankel-norm approximations of linear multivariable
  systems and their $l^{\infty}$-error bounds.
\newblock {\em Internat. J. Control}, 39(6):1115--1193, 1984.

\bibitem{GoyR17}
P.~Goyal and M.~Redmann.
\newblock Towards time-limited $\mathcal{H}_2$-optimal model order reduction,
  2017.

\bibitem{GugA03}
A.C. Gugercin, S.and~Antoulas.
\newblock A time-limited balanced reduction method.
\newblock In {\em Proceedings of the 42nd IEEE Conference on Decision and
  Control}. IEEE, 2003.

\bibitem{gugercin2001approximation}
S.~Gugercin, A.C. Antoulas, and N~Bedrossian.
\newblock Approximation of the international space station 1r and 12a models.
\newblock In {\em Decision and Control, 2001. Proceedings of the 40th IEEE
  Conference on}, volume~2, pages 1515--1516. IEEE, 2001.

\bibitem{GugBA08}
S.~Gugercin, C.~Beattie, and A.~C. Antoulas.
\newblock $\mathcal{H}_2$ model reduction for large-scale linear dynamical
  systems.
\newblock {\em Siam J. Matrix Anal. Appl.}, 30(2):609--638, 2008.

\bibitem{Hal92}
Y.~Halevi.
\newblock Frequency weighted model reduction via optimal projection.
\newblock {\em IEEE Transactions on automatic control}, 37(10):1537--1542,
  1992.

\bibitem{Hess2014}
M.W. Hess and P.~Benner.
\newblock A reduced basis method for microwave semiconductor devices with
  geometric variations.
\newblock {\em COMPEL: The International Journal for Computation and
  Mathematics in Electrical and Electronic Engineering}, 33(4):1071--1081,
  2014.

\bibitem{HolLB96}
P.~Holmes, J.~L. Lumley, and G.~Berkooz.
\newblock {\em Turbulence, coherent structures, dynamical systems and symmetry
  (Cambridge Monographs on Mechanics)}.
\newblock Cambridge University Press, Cambridge, UK, 1996.

\bibitem{HylB85}
D.C. Hyland and D.S. Bernstein.
\newblock The optimal projection equations for model reduction and the
  relationships among the methods of wilson, skelton, and moore.
\newblock {\em IEEE Trans. Automat. Control,.}, 30:1201--1211, 1985.

\bibitem{kellems2009low}
A.R. Kellems, D.~Roos, N.~Xiao, and S.J. Cox.
\newblock Low-dimensional, morphologically accurate models of subthreshold
  membrane potential.
\newblock {\em Journal of computational neuroscience}, 27(2):161, 2009.

\bibitem{Kur17}
P.~K\"urschner.
\newblock Balanced truncation model order reduction in limited time intervals
  for large systems.
\newblock Advances in Computational Mathematics, 2018.

\bibitem{LepMPV91}
A.~Lepschy, G.A. Mian, G.~Pinato, and U.~Viaro.
\newblock Rational $l_2$ approximation: A nongradient algorithm.
\newblock {\em Proceedings of the 30th IEEE Conference on Decision and
  Control}, pages 2321--2323, 1991.

\bibitem{Lieberman2010}
C.~Lieberman, K.~Willcox, and O.~Ghattas.
\newblock Parameter and state model reduction for large-scale statistical
  inverse problems.
\newblock {\em SIAM Journal on Scientific Computing}, 32(5):2523--2542, August
  2010.

\bibitem{magruder2010rational}
C.~Magruder, C.A. Beattie, and S.~Gugercin.
\newblock Rational krylov methods for optimal $\mathcal{L}_2$ model reduction.
\newblock In {\em Decision and Control (CDC), 2010 49th IEEE Conference on},
  pages 6797--6802. IEEE, 2010.

\bibitem{MeiL67}
L.~Meier and D.G. Luenberger.
\newblock Approximation of linear constant systems.
\newblock {\em IEE. Trans. Automat. Contr.}, 12:585--588, 1967.

\bibitem{MelVG14}
S.~A. Melchior, P.~Van~Dooren, and K.~A. Gallivan.
\newblock Model reduction of linear time-varying systems over finite horizons.
\newblock {\em Applied Numerical Mathematics}, 77:72--81, 2014.

\bibitem{Moo81}
B.~Moore.
\newblock Principal component analysis in linear systems: Controllability,
  observability, and model reduction.
\newblock {\em IEEE Transactions on Automatic Control}, 26(1):17--32, 1981.

\bibitem{MulR76}
C.T. Mullis and R.~Roberts.
\newblock Synthesis of minimum roundoff noise fixed point digital filters.
\newblock {\em Circuits and Systems, IEEE Transactions on}, 23(9):551--562,
  1976.

\bibitem{panzer2013}
H.K.F. Panzer, S.~Jaensch, T.~Wolf, and B.~Lohmann.
\newblock A greedy rational {K}rylov method for $\ch_2$-pseudooptimal model
  order reduction with preservation of stability.
\newblock In {\em American Control Conference (ACC), 2013}, pages 5512--5517,
  2013.

\bibitem{RedK17}
P.~Redmann, M.and~K\"urschner.
\newblock An $\mathcal{H}_2$-type error bound for time-limited balanced
  truncation, 2017.

\bibitem{SpaMM92}
J.T. Spanos, M.H. Milman, and D.L. Mingori.
\newblock A new algorithm for $\mathcal{L}_2$ optimal model reduction.
\newblock {\em Automat.}, 28:897--909, 1992.

\bibitem{VanGA08}
P.~van Dooren, K.A. Gallivan, and P.A. Absil.
\newblock $\mathcal{H}_2$-optimal model reduction of mimo systems.
\newblock {\em Applied Mathematics Letters}, 21(12):1267--1273, 2008.

\bibitem{vuillemin2014poles}
P.~Vuillemin, C.~Poussot-Vassal, and D.~Alazard.
\newblock Poles residues descent algorithm for optimal frequency-limited
  $\mathcal{H}_2$ model approximation.
\newblock In {\em Control Conference (ECC), 2014 European}, pages 1080--1085.
  IEEE, 2014.

\bibitem{Wil70}
D.~A. Wilson.
\newblock Optimum solution of model-reduction problem.
\newblock {\em Proceedings of the Institution of Electrical Engineers},
  117(6):1161--1165, 1970.

\bibitem{YanL99}
W.~Y. Yan and J.~Lam.
\newblock An approximate approach to $\mathcal{H}_2$ optimal model reduction.
\newblock {\em IEEE Trans. Automat. Control}, 44:1341--1358, 1999.

\end{thebibliography}
\end{document}